\documentclass{amsart}
\usepackage{amsmath,amsthm,amssymb,IMjournal}
\begin{document}
\newtheorem{The}{Theorem}[section]
\newtheorem{Def}{Definition}[section]
\newtheorem{lem}{Lemma}[section]

\newtheorem{proposition}{Proposition}[section]
\numberwithin{equation}{section}

\title{The First Cohomology of Lie Superalgebra $\tilde{P}(2)$ with Coefficients in Kac Modules}

\author{\|Zilu Zhang$^{1}$ |Liping Sun$^{1}$ |Zhaoxin Li$^{2}$\\1.College of Science, Harbin University of Science and Technology, Harbin 150080, China;
\\
        2.School of Mathematics Sciences, Heilongjiang University, Harbin 150080, China.}

\rec {October 14, 2022}

\dedicatory{Cordially dedicated to ...}

\abstract
%%%%%%
%%%7 to 12 lines optimum; references in the abstract should be given in full form, e.g. V. Nov\'ak, M. Novotn\'y: Linear extensions of orderings. Czech. Math. J. 50 (2000), 853--864.
%%%%%
   Over a field of characteristic $p>2$, firstly, the structure of Kac modules of Lie superalgebra $\tilde{P}(2)$ and the weight space decompositions are given.
  Secondly, the weight-derivations of $\tilde{P}(2)$ to its Kac modules are computed. Finally, the first cohomology of $\tilde{P}(2)$ with coefficients in Kac modules is determined.
\endabstract

\keywords
   Lie superalgebras, Kac modules, derivation, cohomology
\endkeywords

\subjclass
%%%%%
%%%Mathematics Subject Classification 2010
%%%%%
17B05, 17B50
\endsubjclass

\thanks
   Corresponding author(Liping Sun): sunliping@hrbust.edu.cn.\\
   \indent The research has been supported by the National Natural Science Foundation of China [grant numbers 12061029] and the Natural Science Foundation of Heilongjiang Province of China [grant numbers QC2018006].
\endthanks

\section{introduction}\label{sec1}

It is well known that  the representation is of great importance in the study of both Lie algebras and Lie superalgebras. In the present, the representation of Lie algebras has a relatively complete system. Lie superalgebras are the natural generalization of Lie algebras and can be divided into modular Lie superalgebras and non-modular Lie superalgebras according to the different characteristics of basic fields.
Since the 1970s, many important research achievements have been get in the representations of non-modular Lie superalgebras, such as \cite{ref1, ref2, ref3, ref4}. So many researchers began to focus on the representation of modular Lie superalgebras \cite{ref13, ref14, ref15, ref16, ref17, ref18, ref19}.
The cohomology groups are helpful for the representations of modular Lie superalgebras. First Leits and Fuks calculated the cohomology of classical Lie superalgebras with coefficients in trivial modules in 1984. Then  the following  results on the Lie superalgebras cohomology  are gradually plentiful and substantial, include some  discussion on nontrivial modules, for example \cite{ref6, ref7, ref8, ref9}. But very little seems to be known about the study of modular Lie superalgebras cohomology with coefficients in  nontrivial or adjoint modules except \cite{ref9, ref10, ref11, ref12}.

In representation theories of Lie superalgebras, the  Kac modules are of significance. In 2007, Su Yucai and Zhang Ruibin calculated the first and second cohomology of $\mathfrak{sl}_{m|n}$ and $\mathfrak{osp}_{2|2n}$ to finite-dimensional irreducible modules and Kac modules\cite{ref7}. In 2010, Shu Bin and Zhang Chaowen studied Witt superalgebras and defined its Kac modules\cite{ref8}.
In 2020, Wang Shujuan and Liu Wende studied the first cohomology of $\mathfrak{sl}_{2|1}$ with coefficients in $\chi$-reduced Kac modules and simple modules \cite{ref9}. In the present article, over a field of characteristic $p>2$, we describe firstly the structure of Kac modules of Lie superalgebra $\tilde{P}(2)$ and compute the weight space decompositions of $\tilde{P}(2)$ and its  Kac modules relative to a fixed Cartan subalgebra  $\mathfrak{h}$ of $\tilde{P}(2)$. Then the work under consideration is reduced to computing the weight-derivations, which preserve the $\mathfrak{h}^*$-gradings, of $\tilde{P}(2)$ to these weight spaces. Finally, the first cohomology of Lie superalgebra $\tilde{P}(2)$ with coefficients in Kac modules are obtained by means of the fact that each derivation of a finite dimensional Lie superalgebras  to its module is equal to a weight-derivation module an inner derivation.

\section{preliminaries}\label{sec2}
Throughout the paper, all vector spaces are over a field $\mathbb{F}$ of characteristic $p>2$ and
finite-dimensional. Let $\mathbb{Z}_{2}:=\{\bar{0}, \bar{1}\}$ be the two-element field.
The symbol $|x|$ or ${\rm {zd}}$(\emph{x}) denotes the $\mathbb{Z}_{2}$-degree or $\mathbb{Z}$-degree of a $\mathbb{Z}_{2}$-homogeneous or $\mathbb{Z}$-homogeneous element $x$  respectively. And we said that element $x$ is even or odd, if $|x|$ is $\bar{0}$ or $\bar{1}$.
The set of all $\mathbb{Z}_{2}$-homogeneous elements in the $\mathbb{Z}_{2}$-degree space $V$ is represented by hg($V$).
Write $\langle v_{1},\ldots,v_{k}\rangle$ for the vector space spanned by $v_{1},\ldots,v_{k}$ over $\mathbb{F}$. In the $\mathbb{Z}_{2}$-degree vector space
$\langle v_{1},\ldots, v_{m}\mid v_{m+1},\ldots, v_{m+k}\rangle,$  we assume that $|v_{i}|=\bar{0}$ and $|v_{m+j}|=\bar{1}$, where $i=1,\ldots,m,$  $j=1,\ldots,k.$

Let $M$ be a Lie superalgebra $\mathfrak{g}$-module. Recall that a $\mathbb{Z}_{2}$-homogeneous linear mapping $\varphi$ is a derivation of parity $|\varphi|$ of $\mathfrak{g}$ to $\mathfrak{g}$-module $M$ provided that
$$\varphi ([x,y]) = {( - 1)^{\left| \varphi  \right|\left| x \right|}}x \varphi (y) - {( - 1)^{\left| y \right|(\left| \varphi  \right| + \left| x \right|)}}y \varphi (x),\,\textmd{for}\,\textmd{all}\,x,y \in \textmd{hg}(\mathfrak{g}).$$
A derivation $\varphi$ of $\mathfrak{g}$ to $M$ is said to be inner determined by $v$ if there exists $v\in\textmd{hg}(\mathfrak{g})$ such that
$\varphi$(\emph{x})=($-1$)$^{{\mid{x}\mid}{\mid{v}\mid}}$\emph{x}\emph{v} for any $\emph{x}\in\textmd{hg}(\mathfrak{g})$, record as $\mathfrak{D}_{v}$.
Otherwise, $\varphi$ is called an outer derivation. Let $\mathfrak{h}$ be a Cartan subalgebra in the even part of $\mathfrak{g}$. Suppose that $\mathfrak{g}$ and $M$ possess weight space decompositions with respect to $\mathfrak{h}$: $\mathfrak{g}$=$\oplus${$_{\alpha\in\mathfrak{h}^{*}}$}$\mathfrak{g}_{\alpha}$ and
$M$=$\oplus${$_{\alpha\in\mathfrak{h}^{*}}$}\emph{M}$_{\alpha}$.
A derivation $\varphi$ of $\mathfrak{g}$ to $M$ is called a weight-derivation relative to $\mathfrak{h}$ if $\varphi$\emph{($\mathfrak{g}$}$_{\alpha}$)$\subseteq\!M_{\alpha}$, for all {$\alpha${$\in$$\mathfrak{h}$$^{*}$}}.
Let Der($\mathfrak{g},M$) denote the vector space spanned by all the $\mathbb{Z}_{2}$-homogeneous derivations of $\mathfrak{g}$ to $M$. Write Ider($\mathfrak{g},M$) for the vector space spanned by all inner derivations.
The first cohomology of $\mathfrak{g}$ with coefficients in $M$ is the quotient module:
$${\rm H}^{1}(\mathfrak{g},M)=\rm Der(\mathfrak{g},\emph{M})/Ider(\mathfrak{g},\emph{M}). $$

\begin{Def}\label{ddd}\rm\textsuperscript{\cite{ref7}}
Let $\mathfrak{g}=\mathfrak{g}_{-1}\oplus\mathfrak{g}_{0}\oplus\mathfrak{g}_{+1}$ be a restricted Lie superalgebra with $\mathbb{Z}$-grading.
Suppose that $M(\lambda)$ is a simple  finite-dimensional $\mathfrak{g}_{0}$-module with the highest weight $\lambda$, and $\mathfrak{g}_{+1}M(\lambda)=0$. Regarding $M(\lambda)$ as $(\mathfrak{g}_{0}\oplus\mathfrak{g}_{+1})$-module, we call the induced module
$$ K(\lambda)=U(\mathfrak{g})\otimes_{U(\mathfrak{g}_{0}\oplus\mathfrak{g}_{+1})}M(\lambda)$$
restricted $\emph{Kac}$ module of $\mathfrak{g},$ where $U(\mathfrak{g})$  is the enveloping algebra of Lie superalgebra $\mathfrak{g}.$
\end{Def}

Note that $K(\lambda)\cong U(\mathfrak{g}_{-1})\otimes_{\mathbb{F}}M(\lambda)$ as a vector space.

\section{the structure of kac modules}\label{sec3}
Lie superalgebra $\tilde{P}(2)$ is defined as follows:
\[ \tilde P(2): = \left\{ {\left. {\left. {\left( {\left. {\begin{array}{*{20}{c}}
\emph{A}&B\\
C&{ - {A^T}}
\end{array}} \right) \in \emph{gl}(2,2)} \right.} \right|\emph{B} = {\emph{B}^\emph{T}},\emph{C} =  - {\emph{C}^\emph{T}}} \right\}} \right..\]
From now on, write $\mathfrak{g}$ for $\tilde{P}(2)$ and  $e_{ij}$ for the 4 $\times$ 4 matrices  which has 1 in the position $(i,j)$ and 0 elsewhere.
 Let $$\gamma=e_{41}-e_{32},\,\,h_{1}=e_{33}-e_{11},\,\,h_{2}=e_{44}-e_{22},\,\,
\alpha=e_{43}-e_{12},\,\,\beta=e_{34}-e_{21}.$$
Note that $\mathfrak{g}$ possesses a $\mathbb{Z}$-grading structure $\mathfrak{g}=\mathfrak{g}_{-1}\oplus\mathfrak{g}_{0}\oplus\mathfrak{g}_{+1}$, then the following elements form the basis of $\mathfrak{g}$:
$$\mathfrak{g}_{-1}=\langle\gamma\rangle,\,\mathfrak{g}_{0}=\langle h_{1},\,h_{2},\,\alpha,\,\beta\rangle,\,
\mathfrak{g}_{+1}=\langle e_{13},\ e_{24},\ e_{14}+e_{23}\rangle. $$
Fix the standard Cartan subalgebra $\mathfrak{h}$ of $\mathfrak{g}_{\bar{0}}$ spanned by $h_{1}$ and $h_{2}$.
Let $\varepsilon_{1}, \varepsilon_{2}\in\mathfrak{h}^{*}$ make
$\varepsilon_{i}(h_{j})=\delta_{ij},\,i, j=1, 2$. Then the roots and the root-vectors for $\mathfrak{g}$ can be obtained in the Table 3.1 below.
\begin{center}
Table 3.1: Roots and root-vectors for $\mathfrak{g}$
\end{center}
\begin{center}
 \begin{tabular}{c c c c c c c c  }
           \hline

           Roots & $ \theta $& $-2\varepsilon_{1}$ & $-\varepsilon_{1}-\varepsilon_{2}$ & $-\varepsilon_{1}+\varepsilon_{2}$ & $-2\varepsilon_{2}$
           & $\varepsilon_{1}-\varepsilon_{2}$ & $\varepsilon_{1}+\varepsilon_{2}$\\
           \hline
           Root-vectors & $ h_{1}, h_{2}$ & $e_{13}$ &$e_{14}+e_{23}$ & $ \alpha$ & $e_{24}$ & $\beta$ & $\gamma$\\
           \hline
         \end{tabular}
         \end{center}

Let $\lambda=a\varepsilon_{1}+b\varepsilon_{2}$, where $a,b\in\mathbb{F}_{p}$.
Suppose $v_{0}$ satisfies: $h_{i}v_{0}=\lambda(h_{i})v_{0},\,\alpha v_{0}=0$.  Inductively define $v_{k}=\beta^{k}v_{0}$, where $k\in\mathbb{Z}$, and note that $\beta^{[p]}=0$, then $\beta^{p}v_{0}=0$, that is, $v_{p}=0$. Therefore, we get that $v_{0}, v_{1}, \ldots, v_{p-1}$ are nonzero and linearly independent. Then $V:=\langle v_{0}, v_{1}, \ldots, v_{p-1}\rangle$ is a $p$-dimensional module of $\mathfrak{g}_{0}$.
For any $c\in\mathbb{F}_{p}$, let $\Phi(c)\in \{0, 1, \cdots, p-1\}\subseteq \mathbb{Z}$  such that $\Phi(c)\equiv c(\textmd{mod}p).$  Also due to
$$h_{1}v_{k}=(a+k)v_{k},\,h_{2}v_{k}=(b-k)v_{k},$$
$$\alpha v_{k}=k(b-a-k+1)v_{k-1},\,\beta v_{k}=v_{k+1},\,0\leq k\leq p-1.$$
It is easy to know $W:=\langle v_{\Phi(b-a+1)}, \ldots, v_{p-1}\rangle$ is a maximal submodule of $V$.
Hence $M(\lambda):=V/W$ becomes a simple module of $\mathfrak{g}_{0}$.
Without confusion, we write the images of the elements of $V$ in $M(\lambda)$ still by the elements of $V$ itself. Thus
$$M(\lambda)=\langle v_{0}, v_{1}, \ldots, v_{\Phi(b-a)}\rangle.$$
Let $\mathfrak{g}_{+1}M(\lambda)=0$, $M(\lambda)$ be regarded as the simple module of $\mathfrak{g}_{0}\oplus\mathfrak{g}_{+1}.$
By Definition \ref{ddd},
$$K(\lambda)=\langle1\otimes v_{k}, \gamma\otimes v_{k}, k=0, 1, \ldots ,{\Phi(b-a)}\rangle.$$
It is easy to see that $|1\otimes v_{k}|=\bar{0},\,|\gamma\otimes v_{k}|=\bar{1}$. The module action of $\mathfrak{g}$ on $K(\lambda)$ is given below:
$$h_{1}(1\otimes v_{k})=(a+k)\otimes v_{k},\,h_{2}(1\otimes v_{k})=(b-k)\otimes v_{k},$$
$$h_{1}(\gamma\otimes v_{k})=(a+k+1)\otimes v_{k},\,h_{2}(\gamma\otimes v_{k})=(b-k+1)\otimes v_{k},$$
$$\alpha(1\otimes v_{k})=k(b-a-k+1)\otimes v_{k-1},\,\alpha(\gamma\otimes v_{k})=k(b-a-k+1)\gamma\otimes v_{k-1},$$
$$\beta(1\otimes{v_k}) = \;\left\{ {\begin{array}{*{2}{ll}}
{1\otimes v_{k+1},\quad 0 \le k < \Phi(b-a)},\\
{\quad 0,\quad\quad\quad\quad k = \Phi(b-a),}
\end{array}} \right.
\beta(\gamma\otimes{v_k}) = \;\left\{ {\begin{array}{*{2}{ll}}
{\gamma\otimes v_{k+1},\quad 0 \le k < \Phi(b-a)},\\
{\quad 0,\quad\quad\quad\quad k = \Phi(b-a),}
\end{array}} \right.$$

$$ e_{13}(1\otimes v_{k})=e_{24}(1\otimes v_{k})=(e_{14}+e_{23})(1\otimes v_{k})=0,$$
$$\,\gamma(1\otimes v_{k})=\gamma\otimes v_{k},\,\gamma(\gamma\otimes v_{k})=0,$$
$$ e_{13}(\gamma\otimes v_{k})=k(b-a-k+1)\otimes v_{k-1},\,
(e_{14}+e_{23})(\gamma\otimes v_{k})=(b-a-2k)\otimes v_{k},,$$

$$e_{24}(\gamma\otimes{v_k}) = \;\left\{ {\begin{array}{*{2}{ll}}
{-1\otimes v_{k+1},\quad 0 \le k < \Phi(b-a)},\\
{\quad 0,\quad\quad\quad\quad k = \Phi(b-a).}
\end{array}} \right.$$
We list the weight-vector of $K(\lambda)$ relative to $\mathfrak{h}$ , where $a,\,b\in\mathbb{F}_{p}$ (see Table 3.2).
\begin{center}
Table 3.2: Weights and weight-vectors of $K(\lambda)$
\end{center}
\begin{center}
 \begin{tabular}{c c c }
           \hline

           Weights & $(a+k)\varepsilon_{1}+(b-k)\varepsilon_{2}$ & $(a+k+1)\varepsilon_{1}+(b-k+1)\varepsilon_{2}$ \\
           \hline
           Weight-vectors & $1\otimes v_{k}$ & $\gamma\otimes v_{k}$ \\
           \hline
         \end{tabular}
         \end{center}

\section{target-weight spaces of $K(\lambda)$}\label{sec4}
The following fact will simplify the computation of the first cohomology $H^{1}(\mathfrak{g}, K(\lambda))$.

\begin{lem}\label{y3.1}\rm\textsuperscript{\cite{ref9, ref10}}
Each derivation of finite dimensional Lie superalgebra $\mathfrak{g}$ to $\mathfrak{g}$-module M is equal to a weight-derivation module an inner derivation.
 \end{lem}

In view of Lemma \ref{y3.1} and the definition of the first cohomology, it is sufficient to compute the weight-derivations of $\mathfrak{g}$ to $K(\lambda)$ relative to $\mathfrak{h}$, if we want to consider the first cohomology of $\mathfrak{g}$ to $K(\lambda)$.
However, the target-weight spaces should be given before calculating the weight-derivation. This section aims to solve this problem.
Let $c,\,d\in\mathbb{Z}$, $\overline{c, d}:=\{x\in\mathbb{Z}|c\leq x \leq d\}.$

\begin{lem}\label{y3.2}
The relationships about $\Phi(b), \Phi(b+1), \Phi(b+2)$ and $\Phi(2b+2)$ are as follows $:$\\
$(\emph{i})\, \Phi(b)\leq\Phi(2b+2)\Longleftrightarrow \Phi(b)\in \overline{0, \frac{p-3}{2}}\cup\{p-2\}$, \\
$(\emph{ii})\, \Phi(b+1)\leq\Phi(2b+2)\Longleftrightarrow \Phi(b)\in \overline{0, \frac{p-3}{2}}\cup\{p-1\}$, \\
$(\emph{iii})\, \Phi(b+2)\leq\Phi(2b+2)\Longleftrightarrow \Phi(b)\in \overline{0, \frac{p-3}{2}}\cup\{p-2\}$. \\
\end{lem}
\begin{proof}
Since $\Phi(b)\in \overline{0,p-1}$, the scopes of $\Phi(b), \Phi(b+1), \Phi(b+2)$ and $\Phi(2b+2)$ can be obtained(see Table 4.1).
 \begin{center}
Table 4.1: The scopes of $\Phi(b), \Phi(b+1), \Phi(b+2)$ and $\Phi(2b+2)$
\end{center}
\begin{center}
 \begin{tabular}{c c c c c c c c c c}
           \hline

           $\Phi(b)$   &0 & 1 &$\ldots$ &$\frac{p-3}{2}$ &$\frac{p-1}{2}$ &$\ldots$ &$p-3$ &$p-2$ &$p-1$\\
           \hline
           $\Phi(b+1)$ &1 & 2 &$\ldots$ &$\frac{p-1}{2}$ &$\frac{p+1}{2}$ &$\ldots$ &$p-2$ &$p-1$ &$0$\\
           \hline
           $\Phi(b+2)$ &2 & 3 &$\ldots$ &$\frac{p+1}{2}$ &$\frac{p+3}{2}$ &$\ldots$ &$p-1$ &0 &1\\
           \hline
           $\Phi(2b+2)$&2 & 4 &$\ldots$ &$p-1$ &1 &$\ldots$ &$p-4$ &$p-2$ &0\\
           \hline
         \end{tabular}
         \end{center}

It can be seen from Table 4.1 that
$$\Phi(b)\in \overline{0, \frac{p-3}{2}}\cup\{p-2\},\,if\,\Phi(b)\leq\Phi(2b+2),$$
$$\Phi(b)\in \overline{0, \frac{p-3}{2}}\cup\{p-1\},\,if\,\Phi(b+1)\leq\Phi(2b+2),$$
$$\Phi(b)\in \overline{0, \frac{p-3}{2}}\cup\{p-2\},\,if\,\Phi(b+2)\leq\Phi(2b+2).$$\end{proof}

By the similar methods, we can proof the following  Lemma \ref{y3.3}-Lemma \ref{y3.5}.

\begin{lem}\label{y3.3}
The relationships about $\Phi(b), \Phi(b-1), \Phi(b+1)$ and $\Phi(2b)$ are as follows $:$\\
$(\emph{i})\, \Phi(b)\leq\Phi(2b)\Longleftrightarrow \Phi(b)\in \overline{0, \frac{p-1}{2}}$, \\
$(\emph{ii})\, \Phi(b-1)\leq\Phi(2b)\Longleftrightarrow \Phi(b)\in \overline{1, \frac{p-1}{2}}\cup\{p-1\}$, \\
$(\emph{iii})\, \Phi(b+1)\leq\Phi(2b)\Longleftrightarrow \Phi(b)\in \overline{1, \frac{p-1}{2}}\cup\{p-1\}$. \\
\end{lem}

\begin{lem}\label{y3.4}
The relationships about $\Phi(b+1), \Phi(b+2), \Phi(b+3)$ and $\Phi(2b+4)$ are as follows $:$\\
$(\emph{i})\, \Phi(b+1)\leq\Phi(2b+4)\Longleftrightarrow \Phi(b)\in \overline{0, \frac{p-5}{2}}\cup\{p-1, p-3\}$, \\
$(\emph{ii})\, \Phi(b+2)\leq\Phi(2b+4)\Longleftrightarrow \Phi(b)\in \overline{0, \frac{p-5}{2}}\cup\{p-1, p-2\}$, \\
$(\emph{iii})\, \Phi(b+3)\leq\Phi(2b+4)\Longleftrightarrow \Phi(b)\in \overline{0, \frac{p-5}{2}}\cup\{p-1, p-3\}$. \\
\end{lem}

\begin{lem}\label{y3.5}
The relationships about $\Phi(b-1)$ and $\Phi(2b-2)$ are as follows $:$
$\Phi(b-1)\leq\Phi(2b-2)\Longleftrightarrow \Phi(b)\in \overline{1, \frac{p+1}{2}}$.
\end{lem}
The above work is all in preparation for finding the target-weight spaces of $K(\lambda)$.

\begin{proposition}\label{m3.1}
The target-weight spaces of $K(\lambda)$ are
\begin{equation}\label{equ1}
K(\lambda)_{(-2, 0)}=\left\{
\begin{aligned}
&\langle1\otimes v_{\Phi(b)}\mid0 \rangle,\,\Phi(b)\in \overline{0, \frac{p-3}{2}}\cup\{p-2\}\,and\,a+b=-2,\\
&\langle0\mid\gamma\otimes v_{\Phi(b+1)}\rangle,\,\Phi(b)\in \overline{0, \frac{p-5}{2}}\cup\{p-1, p-3\}\,and\,a+b=-4,\\
&\langle0\mid0\rangle,\,otherwise.
\end{aligned}
\right.
\end{equation}

\begin{equation}\label{equ2}
K(\lambda)_{(-1, -1)}=\left\{
\begin{aligned}
&\langle1\otimes v_{\Phi(b+1)}\mid0\rangle,\,\Phi(b)\in \overline{0, \frac{p-3}{2}}\cup\{p-1\}\,and\,a+b=-2,\\
&\langle0\mid\gamma\otimes v_{\Phi(b+2)}\rangle,\,\Phi(b)\in \overline{0, \frac{p-5}{2}}\cup\{p-1, p-2\}\,and\,a+b=-4,\\
&\langle0\mid0\rangle,\,otherwise.
\end{aligned}
\right.
\end{equation}

\begin{equation}\label{equ3}
K(\lambda)_{(-1, 1)}=\left\{
\begin{aligned}
&\langle1\otimes v_{\Phi(b-1)}\,|\,\, 0\,\, \rangle,\,\Phi(b)\in \overline{1, \frac{p-1}{2}}\cup\{p-1\}\,and\,a+b=0,\\
&\langle0\mid\gamma\otimes v_{\Phi(b)}\rangle,\,\Phi(b)\in \overline{0, \frac{p-3}{2}}\cup\{p-2\}\,and\,a+b=-2,\\
&\langle0\mid0\rangle,\,otherwise.
\end{aligned}
\right.
\end{equation}

\begin{equation}\label{equ4}
K(\lambda)_{(0, -2)}=\left\{
\begin{aligned}
&\langle1\otimes v_{\Phi(b+2)}\mid0\rangle,\,\Phi(b)\in \overline{0, \frac{p-3}{2}}\cup\{p-2\}\,and\,a+b=-2,\\
&\langle0\mid\gamma\otimes v_{\Phi(b+3)}\rangle,\,\Phi(b)\in \overline{0, \frac{p-5}{2}}\cup\{p-1, p-3\}\,and\,a+b=-4,\\
&\langle0\mid0\rangle,\,otherwise.
\end{aligned}
\right.
\end{equation}

\begin{equation}\label{equ5}
K(\lambda)_{(1, -1)}=\left\{
\begin{aligned}
&\langle1\otimes v_{\Phi(b+1)}\mid0\rangle,\,\Phi(b)\in \overline{1, \frac{p-1}{2}}\cup\{p-1\}\,and\,a+b=0,\\
&\langle0\mid\gamma\otimes v_{\Phi(b+2)}\rangle,\,\Phi(b)\in \overline{0, \frac{p-3}{2}}\cup\{p-2\}\,and\,a+b=-2,\\
&\langle0\mid0\rangle,\,otherwise.
\end{aligned}
\right.
\end{equation}

\begin{equation}\label{equ6}
K(\lambda)_{(1, 1)}=\left\{
\begin{aligned}
&\langle1\otimes v_{\Phi(b-1)}\mid0\rangle,\,\Phi(b)\in \overline{1, \frac{p+1}{2}}\,and\,a+b=2,\\
&\langle0\mid\gamma\otimes v_{\Phi(b)}\rangle,\,\Phi(b)\in \overline{0, \frac{p-1}{2}}\,and\,a+b=0,\\
&\langle0\mid0\rangle,\,otherwise.
\end{aligned}
\right.
\end{equation}

\begin{equation}\label{equ7}
K(\lambda)_{(0, 0)}=\left\{
\begin{aligned}
&\langle1\otimes v_{\Phi(b)}\mid0\rangle,\,\Phi(b)\in \overline{0, \frac{p-1}{2}}\,and\,a+b=0,\\
&\langle0\mid\gamma\otimes v_{\Phi(b+1)}\rangle,\,\Phi(b)\in \overline{0, \frac{p-3}{2}}\cup\{p-1\}\,and\,a+b=-2,\\
&\langle0\mid0\rangle,\,otherwise.
\end{aligned}
\right.
\end{equation}
\end{proposition}

\begin{proof}
We take (\ref{equ6}) and (\ref{equ7}) as examples to prove Proposition \ref{m3.1}. Equations (\ref{equ1})-(\ref{equ5}) can be proved similarly.
By Table 2, the weights of $1\otimes v_{k}$ and $\gamma\otimes v_{k}$ are $(a+k,b-k)$ and $(a+k+1,b-k+1)$ respectively.
In order to prove  (\ref{equ6}), we need the weights of $1\otimes v_{k}$ and $\gamma\otimes v_{k}$ to be $(1, 1)$. There are the following two cases.

Case 1: Let $(a+k, b-k)=(1, 1)$. We have $a+b=2$, so that $b-a=2b-2$, thus $\Phi{(b-a)}=\Phi(2b-2)$. The following conclusions can be drawn:\\
Subcase 1.1: $a+b\neq2$. $1\otimes v_{k}$'s weight is not $(1, 1)$, $0\leq k\leq \Phi{(b-a)}$.\\
Subcase 1.2: $a+b=2$. By Lemma \ref{y3.5}, we have:\\
Subcase 1.2.1: When $\Phi{(b)}\in \overline{\frac{p+3}{2}, p-1}\cup\{0\}$, $1\otimes v_{k}$'s weight is not $(1, 1)$, $0\leq k\leq \Phi{(b-a)}$;\\
Subcase 1.2.2: When $\Phi{(b)}\in \overline{1, \frac{p+1}{2}}$, $1\otimes v_{k}$'s weight is $(1, 1)$.

Case 2: Let $(a+k+1, b-k+1)=(1, 1)$. We have $a+b=0$, so that $b-a=2b$, thus $\Phi{(b-a)}=\Phi(2b)$. The following conclusions can be drawn:\\
Subcase 2.1: $a+b\neq0$. $\gamma\otimes v_{k}$'s weight is not $(1, 1)$, $0\leq k\leq \Phi{(b-a)}$.\\
Subcase 2.2: $a+b=0$. By Lemma \ref{y3.3} (i), we have:\\
Subcase 2.2.1: When $\Phi{(b)}\in \overline{\frac{p+1}{2}, p-1}$, $\gamma\otimes v_{k}$'s weight is not $(1, 1)$, $0\leq k\leq \Phi{(b-a)}$;\\
Subcase 2.2.2: When $\Phi{(b)}\in \overline{0, \frac{p-1}{2}}$, $\gamma\otimes v_{k}$'s weight is $(1, 1)$.

Analogously, in order to prove (\ref{equ7}), we need the weights of $1\otimes v_{k}$ and $\gamma\otimes v_{k}$ to be $(0, 0)$. There are the following two cases.

Case 1: Let $(a+k, b-k)=(0, 0)$. We have $a+b=0$, so that $b-a=2b$, thus $\Phi{(b-a)}=\Phi(2b)$. The following conclusions can be drawn:\\
Subcase 3.1: $a+b\neq0$. $1\otimes v_{k}$'s weight is not $(0, 0)$, $0\leq k\leq \Phi{(b-a)}$. \\
Subcase 3.2: $a+b=0$. By Lemma \ref{y3.3} (i), we have:\\
Subcase 3.2.1: When $\Phi{(b)}\in \overline{\frac{p+1}{2}, p-1}$, $1\otimes v_{k}$'s weight is not $(0, 0)$, $0\leq k\leq \Phi{(b-a)}$;\\
Subcase 3.2.2: When $\Phi{(b)}\in \overline{0, \frac{p-1}{2}}$, $1\otimes v_{k}$'s weight is $(0, 0)$.

Case 2: Let $(a+k+1, b-k+1)=(0, 0)$. We have $a+b=-2$, so that $b-a=2b+2$, thus $\Phi{(b-a)}=\Phi(2b+2)$. The following conclusions can be drawn:\\
Subcase 4.1: $a+b\neq-2$. $\gamma\otimes v_{k}$'s weight is not $(0, 0)$, $0\leq k\leq \Phi{(b-a)}$. \\
Subcase 4.2: $a+b=-2$. By Lemma \ref{y3.2} (ii), we have:\\
Subcase 4.2.1: When $\Phi{(b)}\in \overline{\frac{p-1}{2}, p-2}$, $\gamma\otimes v_{k}$'s weight is not $(0, 0)$, $0\leq k\leq \Phi{(b-a)}$;\\
Subcase 4.2.2: When $\Phi{(b)}\in \overline{0, \frac{p-3}{2}}\cup\{p-1\}$, $\gamma\otimes v_{k}$'s weight is $(0, 0)$.

In summary, the proofs of (\ref{equ6}) and (\ref{equ7}) are completed.
\end{proof}
\section{$\rm{H}^{1}(\mathfrak{g}, K(\lambda))$}\label{sec5}
 The following lemma will  simplify the calculation of the first cohomology.
\begin{lem}\label{y4.1}
Suppose that $\phi$ is a weight-derivation of $\mathfrak{g}$ to $K(\lambda)$, we have $x\phi(h_{i})=0,\,i=1, 2$, for all $x\in \mathfrak{g}$.
\end{lem}

\begin{proof}
 We assume that $x\in \mathfrak{g}_{\alpha}$, where $\alpha\in \mathfrak{h}^{*}$, then  $h_{i}x=\alpha(h_{i})x$, $i=1, 2$.
Since $\phi$ is a weight-derivation, we get $h_{i}\phi(x)=\alpha(h_{i})\phi(x)$. By the definition of derivation, the following formula holds:
$$\begin{aligned}
\alpha(h_{i})\phi(x)= \phi(\alpha(h_{i})x) &= \phi([h_{i},x])=h_{i}\phi(x)-(-1)^{|x||\phi|}x\phi(h_{i})\\
 &= \alpha(h_{i})\phi(x)-(-1)^{|x||\phi|}x\phi(h_{i}),
\end{aligned}$$
Therefore, $x\phi(h_{i})=0,\,i=1,2$.
\end{proof}

Notice that the weight-mappings from $\mathfrak{g}$ to $K(\lambda)$ are zero mapping for $\Phi(a+b)\notin \{0, 2, p-2, p-4\}$ from Proposition \ref{m3.1}.

Before computing the first cohomology of $\mathfrak{g}$ with coefficients in $K(\lambda)$, we first introduce four outer derivations. Consider the linear mappings of $\mathfrak{g}$ to $K(\lambda)$.

If $a+b=-2$ and $\Phi(b)=p-2$, we define ${\varphi _1},\,{\varphi_2}$, such that
$$\begin{aligned}
{\varphi _1}&:\, \alpha \mapsto \gamma\otimes v_{p-2},\quad e_{13} \mapsto -1\otimes v_{p-2};\\
{\varphi_2}&:\, \beta \mapsto \gamma\otimes v_{0},\quad e_{24} \mapsto 1\otimes v_{0}.\\
\end{aligned}$$

If $a+b=-2$ and $\Phi(b)=p-1$, we define ${\varphi_3}$, such that
$${\varphi_3}:\, h_{i} \mapsto \gamma\otimes v_{0},\,i=1,2.$$

If $a+b=-4$ and $\Phi(b)=p-1$, we define ${\varphi_4}$, such that
$$\varphi_{4}:\, e_{13}\mapsto 2\gamma\otimes v_{0},\quad e_{24}\mapsto \gamma\otimes v_{2},\quad e_{14}+e_{23}\mapsto -2\gamma\otimes v_{1}.$$
Here we take the convention that, the element of ${\rm Hom}(\mathfrak{g},K(\lambda))$ vanishes on the standard
basis elements of $\mathfrak{g}$ which do not appear. For example ${\varphi _1}(h_{1})=0$, the same below.

\begin{lem}\label{y4.2}
Each ${\varphi _k}$ is both an outer derivation and a weight-derivation for $k =1, 2, 3, 4.$
\end{lem}

\begin{proof}
By the definition of derivation and Proposition \ref{m3.1}, it is obvious that ${\varphi _k}$ is a derivation and weight-derivation for $k =1, 2, 3, 4$. Suppose conversely ${\varphi_k}$ is a nonzero inner derivation given by $v\in K(\lambda)$. By the definition of weight-derivation, the weight of $v$ is (0,0).
For ${\varphi_1},\,{\varphi_2},\,{\varphi_4}$, we know $v=0$ by Proposition \ref{m3.1} (1)-(5), contradictorily. Hence ${\varphi _1},\,{\varphi_2},\,\varphi_{4}$ are outer derivations. For ${\varphi_3}$, we may assume $v=e\gamma\otimes v_{0}$ by Proposition \ref{m3.1} (7), where $e\in\mathbb{F}$. According to the definition of inner derivation,
$$\mathfrak{D}_{v}(h_{1})=h_{1}(e\gamma\otimes v_{0})=0\neq\gamma\otimes v_{0}.$$
Contradictorily. So ${\varphi_3}$ is outer.
\end{proof}

Below, we compute ${\textmd{H}^1}(\mathfrak{g},K(\lambda))$. By Lemma \ref{y3.1}, we only need to compute the weight-derivations of $\mathfrak{g}$ to $K(\lambda)$.
\begin{proposition}\label{d4.1}
\begin{equation}
\emph{H}^{1}(\mathfrak{g}, K(\lambda))=\left\{
\begin{aligned}
&\mathbb{F}\varphi_{1}+\mathbb{F}\varphi_{2},\,a+b=-2\,and\,\Phi(b)=p-2,\\
&\mathbb{F}\varphi_{3},\,a+b=-2\,and\,\Phi(b)=p-1,\\
&\mathbb{F}\varphi_{4},\,a+b=-4\,and\,\Phi(b)=p-1,\\
&0,\,otherwise.
\end{aligned}
\right.
\nonumber
\end{equation}
\end{proposition}

\begin{proof}
According to the range of $a+b$ and Proposition \ref{m3.1}, we proof Proposition \ref{d4.1} for the following four cases. Note in advance that the coefficients $m_{i}$ set below are in $\mathbb{F}$, $i=1,\cdots,7$. Let $\varphi$ be a weight-derivation of $\mathfrak{g}$ to $K(\lambda)$ in each of the following cases.

%1.1
Case 1: $a+b=-2$. The following conclusions can be drawn.

Subcase 1.1: $\Phi(b)\in \overline{0, \frac{p-3}{2}}$. By Proposition \ref{m3.1} (1)-(5) and (7), we may assume $\varphi$:
$$\begin{aligned}
&h_{1}\mapsto m_{1}\gamma\otimes v_{\Phi(b+1)},\\
&h_{2}\mapsto m_{2}\gamma\otimes v_{\Phi(b+1)},\\
&\alpha\mapsto m_{3}\gamma\otimes v_{\Phi(b)},\\
&\beta\mapsto m_{4}\gamma\otimes v_{\Phi(b+2)},\\
&e_{13}\mapsto m_{5}\otimes v_{\Phi(b)},\\
&e_{24}\mapsto m_{6}\otimes v_{\Phi(b+2)},\\
&e_{14}+e_{23} \mapsto m_{7}\otimes v_{\Phi(b+2)}.
\end{aligned}$$
Obviously, $|\varphi|=\bar{1}$. Thus, from the Lemma \ref{y4.1} and the definition of derivation, $\varphi$ is a weight-derivation and the following equations hold:
\begin{equation}
\left\{
\begin{aligned}
&\alpha\varphi(h_{i})=m_{i}\alpha(\gamma\otimes v_{\Phi(b+1)})=0,\,i=1,2,\\
&\varphi([\alpha, \beta])=\alpha\varphi(\beta)-\beta\varphi(\alpha),\\
&\varphi([e_{13}, \gamma])=-e_{13}\varphi(\gamma)-\gamma\varphi(e_{13}),\\
&\varphi([e_{24}, \gamma])=-e_{24}\varphi(\gamma)-\gamma\varphi(e_{24}),\\
&\varphi([e_{14}+e_{23}, \gamma])=-(e_{14}+e_{23})\varphi(\gamma)-\gamma\varphi(e_{14}+e_{23}).\\
\end{aligned}
\nonumber
\right.
\end{equation}
That is,
\begin{equation}
\left\{
\begin{aligned}
&m_{i}(b+1)(b+2)\gamma\otimes v_{\Phi(b)}=0,\,i=1,2,\\
&[m_{4}(b+2)(b+1)-m_{3}]\gamma\otimes v_{\Phi(b+1)}=0,\\
&m_{3}\gamma\otimes v_{\Phi(b)}=-m_{5}\gamma\otimes v_{\Phi(b)},\\
&-m_{4}\gamma\otimes v_{\Phi(b+2)}=-m_{6}\gamma\otimes v_{\Phi(b+2)},\\
&-m_{7}\gamma\otimes v_{\Phi(b+1)}=0.\\
\end{aligned}
\nonumber
\right.
\end{equation}
By solving above equations, we have
\begin{equation}
\left\{
\begin{aligned}
&m_{1}=m_{2}=0,\\
&m_{3}=m_{4}(b+1)(b+2),\\
&m_{5}=-m_{4}(b+1)(b+2),\\
&m_{6}=m_{4},\\
&m_{7}=0.\\
\end{aligned}
\nonumber
\right.
\end{equation}
It is easy to verify that $\varphi=m_{4}\mathfrak{D}_{\gamma\otimes v_{\Phi(b+1)}}$, that is, $\varphi$ is inner.
%1.2

Subcase 1.2: $\Phi(b)=p-2$. By Proposition \ref{m3.1} (1) and (3)-(5), we may assume $\varphi$:
$$\begin{aligned}
&\alpha\mapsto m_{3}\gamma\otimes v_{p-2},\\
&\beta\mapsto m_{4}\gamma\otimes v_{0},\\
&e_{13}\mapsto m_{5}\otimes v_{p-2},\\
&e_{24}\mapsto m_{6}\otimes v_{0}.
\end{aligned}$$Obviously, $|\varphi|=\bar{1}$. Thus, it can be known by the definition of derivation that $\varphi$ is a weight-derivation and the following equations hold:
\begin{equation}
\left\{
\begin{aligned}
&\varphi([e_{13}, \gamma])=-e_{13}\varphi(\gamma)-\gamma\psi(e_{13}),\\
&\varphi([e_{24}, \gamma])=-e_{24}\varphi(\gamma)-\gamma\psi(e_{24}).
\end{aligned}
\nonumber
\right.
\end{equation}
It shows that
\begin{equation}
\left\{
\begin{aligned}
&m_{3}\gamma\otimes v_{p-2}=-m_{5}\gamma\otimes v_{p-2},\\
&-m_{4}\gamma\otimes v_{0}=-m_{6}\gamma\otimes v_{0}.
\end{aligned}
\nonumber
\right.
\end{equation}
We get
\begin{equation}
\left\{
\begin{aligned}
&m_{5}=-m_{3},\\
&m_{6}=m_{4}.
\end{aligned}
\nonumber
\right.
\end{equation}From Lemma \ref{y4.2}, it is easy to see that $\varphi=m_{3}\varphi_{1}+m_{4}\varphi_{2}$, and $\varphi_{1},\,\varphi_{2}$ are outer derivations.
%1.3

Subcase 1.3: $\Phi(b)=p-1$. By Proposition \ref{m3.1} (2) and (7), we may suppose $\varphi$:
$$\begin{aligned}
&h_{1}\mapsto m_{1}\gamma\otimes v_{0},\\
&h_{2}\mapsto m_{2}\gamma\otimes v_{0},\\
&e_{14}+e_{23}\mapsto m_{7}\otimes v_{0}.
\end{aligned}$$Obviously, $|\varphi|=\bar{1}$. According to the definition of derivation, we have the following equations:
\begin{equation}
\left\{
\begin{aligned}
&\varphi([\alpha, \beta])=0,\\
&\varphi([\alpha, e_{24}])=0.
\end{aligned}
\nonumber
\right.
\end{equation}
We obtain
\begin{equation}
\left\{
\begin{aligned}
&(m_{2}-m_{1})\gamma\otimes v_{0}=0,\\
&-m_{7}\otimes v_{0}=0.
\end{aligned}
\nonumber
\quad\Rightarrow\quad
\right.
\left\{
\begin{aligned}
&m_{1}=m_{2},\\
&m_{7}=0.
\end{aligned}
\nonumber
\right.
\end{equation}
Write $m_{1}$, $m_{2}$ for $m,\,m\in\mathbb{F}$. From Lemma \ref{y4.2}, it is easy to know that ${\varphi}=m\varphi_{3}$, and $\varphi$ is an outer derivation when $m\neq0$.
%2.1

Case 2: $a+b=-4$. The following conclusions can be drawn.

Subcase 2.1: $\Phi(b)\in \overline{0, \frac{p-5}{2}}$. By Proposition \ref{m3.1} (1), (2) and (4), we may suppose $\varphi$:
$$\begin{aligned}
&e_{13}\mapsto m_{5}\gamma\otimes v_{\Phi(b+1)},\\
&e_{24}\mapsto m_{6}\gamma\otimes v_{\Phi(b+3)},\\
&e_{14}+e_{23}\mapsto m_{7}\gamma\otimes v_{\Phi(b+2)}.
\end{aligned}$$Obviously, $|\varphi|=\bar{0}$.
The following equations hold from the definition of derivation,
\begin{equation}
\left\{
\begin{aligned}
&\varphi([\alpha, e_{13}])=\alpha\varphi(e_{13}),\\
&\varphi([\alpha, e_{24}])=\alpha\varphi(e_{24}),\\
&\varphi([\alpha, e_{14}+e_{23}])=\alpha\varphi(e_{14}+e_{23}).
\end{aligned}
\nonumber
\right.
\end{equation}
Then,
\begin{equation}
\left\{
\begin{aligned}
&m_{5}(b+1)(b+4)\gamma\otimes v_{\Phi(b)}=0,\\
&-m_{7}\gamma\otimes v_{\Phi(b+2)}=m_{6}(b+2)(b+3)\gamma\otimes v_{\Phi(b+2)},\\
&-2m_{5}\gamma\otimes v_{\Phi(b+1)}=m_{7}(b+2)(b+3)\gamma\otimes v_{\Phi(b+1)}.
\end{aligned}
\nonumber
\quad\Rightarrow\quad
\right.
\left\{
\begin{aligned}
&m_{5}=0,\\
&m_{6}=0,\\
&m_{7}=0.
\end{aligned}
\nonumber
\right.
\end{equation}
Therefore, $\varphi=0$.

%2.2
Subcase 2.2: $\Phi(b)=p-1$. By Proposition \ref{m3.1} (1), (2) and (4), we may assume $\varphi$:
$$\begin{aligned}
&e_{13}\mapsto m_{5}\gamma\otimes v_{0},\\
&e_{24}\mapsto m_{6}\gamma\otimes v_{2},\\
&e_{14}+e_{23}\mapsto m_{7}\gamma\otimes v_{1}.
\end{aligned}$$Obviously, $|\varphi|=\bar{0}$.
These equations are obtained by the definition of derivation,
\begin{equation}
\left\{
\begin{aligned}
&\varphi([\alpha, e_{24}])=\alpha\varphi(e_{24}),\\
&\varphi([\alpha, e_{14}+e_{23}])=\alpha\varphi(e_{14}+e_{23}).
\end{aligned}
\nonumber
\right.
\end{equation}
It follows that
\begin{equation}
\left\{
\begin{aligned}
&-m_{7}\gamma\otimes v_{1}=2m_{6}\gamma\otimes v_{1},\\
&-2m_{5}\gamma\otimes v_{0}=2m_{7}\gamma\otimes v_{0}.
\end{aligned}
\nonumber
\quad\Rightarrow\quad
\right.
\left\{
\begin{aligned}
&m_{5}=2m_{6},\\
&m_{7}=-2m_{6}.
\end{aligned}
\nonumber
\right.
\end{equation}
From Lemma \ref{y4.2}, it is easy to see that ${\varphi}=m_{6}\varphi_{4}$, and ${\varphi}$ is an outer derivation when $m_{6}\neq0$.

%2.3
Subcase 2.3: $\Phi(b)=p-3$. By Proposition \ref{m3.1} (1) and (4), we may suppose $\varphi$:
$$\begin{aligned}
&e_{13}\mapsto m_{5}\gamma\otimes v_{p-2},\\
&e_{24}\mapsto m_{6}\gamma\otimes v_{0}.
\end{aligned}$$Obviously, $|\varphi|=\bar{0}$.
We obtain the following equations by definition of derivation:
\begin{equation}
\left\{
\begin{aligned}
&\varphi([\alpha, e_{13}])=\alpha\varphi(e_{13}),\\
&\varphi([\beta, e_{24}])=\beta\varphi(e_{24}).
\end{aligned}
\nonumber
\right.
\end{equation}
It implies that
\begin{equation}
\left\{
\begin{aligned}
&m_{5}(p-2)\gamma\otimes v_{p-3}=0,\\
&m_{6}\gamma\otimes v_{1}=0.
\end{aligned}
\nonumber
\quad\Rightarrow\quad
\right.
\left\{
\begin{aligned}
&m_{5}=0,\\
&m_{6}=0.
\end{aligned}
\nonumber
\right.
\end{equation}
Consequently, $\varphi=0$.

%2.4
Subcase 2.4: $\Phi(b)=p-2$. By Proposition \ref{m3.1} (2), we may assume $\varphi$:
$$\begin{aligned}
&e_{14}+e_{23}\mapsto m_{7}\gamma\otimes v_{0}.
\end{aligned}$$Obviously, $|\varphi|=\bar{0}$. We have
$$-m_{7}\gamma\otimes v_{0}=\varphi([\beta,e_{13}])=0.$$
Comparing the coefficients gives $m_{7}=0$, so $\varphi=0$.

%3.1
Case 3: $a+b=0$. The following conclusions can be drawn.

Subcase 3.1: $\Phi(b)\in \overline{1, \frac{p-1}{2}}$. By Proposition \ref{m3.1} (3) and (5)-(7), we may assume $\varphi$:
$$\begin{aligned}
&h_{1}\mapsto m_{1}\otimes v_{\Phi(b)},\\
&h_{2}\mapsto m_{2}\otimes v_{\Phi(b)},\\
&\alpha\mapsto m_{3}\otimes v_{\Phi(b-1)},\\
&\beta\mapsto m_{4}\otimes v_{\Phi(b+1)},\\
&\gamma\mapsto m_{8}\gamma\otimes v_{\Phi(b)}.
\end{aligned}$$Obviously, $|\varphi|=\bar{0}$. According to Lemma \ref{y4.1} and the definition of derivation we have the following equations:
\begin{equation}
\left\{
\begin{aligned}
&\alpha\varphi(h_{i})=m_{i}\alpha(1\otimes v_{\Phi(b)})=0,\,i=1,2,\\
&\varphi([e_{13}, \gamma])=e_{13}\varphi(\gamma),\\
&\varphi([e_{24}, \gamma])=e_{24}\varphi(\gamma).
\end{aligned}
\nonumber
\right.
\end{equation}
Then,
\begin{equation}
\left\{
\begin{aligned}
&m_{i}b(b+1)\otimes v_{\Phi(b-1)}=0,\,i=1,2,\\
&m_{3}\otimes v_{\Phi(b-1)}=m_{8}b(b+1)\otimes v_{\Phi(b-1)},\\
&-a_{4}\otimes v_{\Phi(b+1)}=-m_{8}\otimes v_{\Phi(b+1)}.
\end{aligned}
\nonumber
\quad\Rightarrow\quad
\right.
\left\{
\begin{aligned}
&m_{1}=m_{2}=0,\\
&m_{3}=b(b+1)m_{8},\\
&m_{4}=m_{8}.
\end{aligned}
\nonumber
\right.
\end{equation}Easy to verify $\varphi=m_{8}\mathfrak{D}_{1\otimes v_{\Phi(b)}}$, that is, $\varphi$ is inner.

%3.2
Subcase 3.2: $\Phi(b)=0$. By Proposition \ref{m3.1} (6) and (7), we may suppose $\varphi$:
$$\begin{aligned}
&h_{1}\mapsto m_{1}\otimes v_{0},\\
&h_{2}\mapsto m_{2}\otimes v_{0},\\
&\gamma\mapsto m_{8}\gamma\otimes v_{0}.
\end{aligned}$$Obviously, $|\varphi|=\bar{0}$. Thus the following equation holds from the Lemma \ref{y4.1}:
$$\gamma\varphi(h_{i})=m_{i}\gamma\otimes v_{0}=0,\,i=1,\,2.$$
Comparing the coefficients gives $m_{1}=m_{2}=0$. By calculation we find that $m_{8}$ is arbitrary, so it is easy to verify $\varphi=m_{8}\mathfrak{D}_{1\otimes v_{0}}$, that is, $\varphi$ is inner.

%3.3
Subcase 3.3: $\Phi(b)=p-1$. By Proposition \ref{m3.1} (3) and (5), we may suppose $\varphi$:
$$\begin{aligned}
&\alpha\mapsto m_{3}\otimes v_{p-2},\\
&\beta\mapsto m_{4}\otimes v_{0}.
\end{aligned}$$Obviously, $|\varphi|=\bar{0}$. We obtain the following equations by definition of derivation:
$$m_{3}\otimes v_{p-2}=\varphi([e_{13},\gamma])=0,$$
$$-m_{4}\otimes v_{0}=\varphi([e_{24},\gamma])=0.$$
Hence $m_{3}=m_{4}=0$, and $\varphi=0$.

%11
Case 4: $a+b=2$ and $\Phi(b)\in \overline{1. \frac{p+1}{2}}$. By Proposition \ref{m3.1} (6), we may suppose $\varphi$:
$$\begin{aligned}
&\gamma\mapsto m_{8}\otimes v_{\Phi(b-1)}.
\end{aligned}$$Obviously, $|\varphi|=\bar{1}$. We get the equation
$$0=\varphi([\gamma,\gamma])=-2\varphi(\gamma)=-2m_{8}\gamma(1\otimes v_{\Phi(b-1)})=-2m_{8}\gamma\otimes v_{\Phi(b-1)}.$$
Therefore, $m_{8}=0$, and $\varphi=0$.

In summary, we get this proposition by the definition of the first cohomology.
\end{proof}

\begin{The}\label{d1.1}
\begin{equation}
\emph{dim}(\emph{H}^{1}(\mathfrak{g}, K(\lambda)))=\left\{
\begin{aligned}
&2,\,if\,a+b=-2 \,and\,\Phi(b)=p-2,\\
&1,\,if\,a+b=-2\,or-4\,and\,\Phi(b)=p-1,\\
&0,\,otherwise,
\end{aligned}
\right.
\nonumber
\end{equation}where $a,b\in\mathbb{F}_{p}$, $\Phi(b)\in \{0, 1, \ldots, p-1\}$.
\end{The}
\begin{proof}
By Proposition \ref{d4.1}, it suffices to prove $\varphi_{1}$ and $\varphi_{2}$ are linearly independent modulo the inner derivation space $\rm{Ider}(\mathfrak{g}, K(\lambda))$  when  $a+b=-2$ and $\Phi(b)=p-2$. Suppose that $0\neq t\in K(\lambda)$, such that $\varphi_{1}$ and $\varphi_{2}$ are linearly dependent modulo $\rm{Ider}(\mathfrak{g}, K(\lambda)).$ Because $\varphi_{1}$ and $\varphi_{2}$ are weight-derivations, then $t\in K(\lambda)_{(0,0)}$. By Proposition \ref{m3.1} (7), we know $t=0$, contradictorily. Other cases, it is obvious  from Proposition \ref{d4.1}.
\end{proof}


{\small
\begin{thebibliography}{999}
%%%%%
%%%References should be listed in alphabetical order; please supply the zbMATH and MathSciNet numbers as indicated if possible.
%%%%%
\bibitem{ref18} F.F. Duan, B. Shu, Y.F. Yao:
 On blocks in restricted representations of Lie superalgebras of Cartan type. {\it Acta. Math. Sin.-English Ser.} {\bf36} (2020) 1061-1075. Zbl 1486.17013, MR4145698, DOI 10.1007/s10114-020-9469-z

\bibitem{ref1}{\it C. Gruson, V. Serganova}:
 Cohomology of generalized supergrassmannians and character formulae for basic classical Lie superalgebras. Proc. Lond. Math. Soc. {\bf101} (2010), 852-892. Zbl 1216.17005, MR2734963, DOI 10.1112/plms/pdq014

\bibitem{ref2}{\it I.M. Musson, V. Serganova}:
 Combinatorics of character formulas for the Lie superalgebra $\mathfrak{gl}(m, n)$. Transform. Groups. {\bf16} (2011) 555-578. Zbl 1275.17018, MR2806501, DOI 10.1007/s00031-011-9147-4

\bibitem{ref19}{\it L. Pan, B. Shu}:
 Jantzen filtration and strong linkage principle for modular Lie superalgebras. Forum Math. {\bf30} (2018) 1573-1598. Zbl 1404.17027, MR3871461, DOI 10.1515/forum-2018-0065

\bibitem{ref4}{\it Y.C. Su}:
 Composition factors of Kac modules for the general linear Lie superalgebras. Math. Z. {\bf252} (2006), 731-754.  Zbl 1106.17008, MR2206623, DOI 10.1007/s00209-005-0874-x

\bibitem{ref11}{\it L.P. Sun, W.D. Liu}:
 Low-dimensional cohomology of Lie superalgebra $A(1, 0)$ with coefficients in Witt or special superalgebras. Taiwanese J. Math. {\bf17} (2013) 83-107. Zbl 1280.17021, MR3028859, DOI 10.11650/tjm.17.2013.1812

\bibitem{ref10}{\it L.P. Sun, W.D. Liu, B.Y. Wu}:
 Low-dimensional cohomology of Lie superalgebra $\mathfrak{sl}_{m\mid n}$ with coefficients in Witt or special superalgebras. Indag. Math. {\bf25} (2014) 59-77. Zbl 1382.17012, MR3131765, DOI 10.1016/j.indag.2013.07.006

\bibitem{ref6}{\it M. Scheunert, R.B. Zhang}:
 Cohomology of Lie Superalgebras and Their Generalizations. J. Math. Phys. {\bf39} (1998) 5024–5061. Zbl 0928.17023, MR1643330,
DOI 10.1063/1.532508

\bibitem{ref3}{\it Y.C. Su, R.B. Zhang}:
 Character and dimension formulae for general linear superalgebra. Adv. Math. {\bf211} (2007), 1-33.
Zbl 1166.17002, MR2313526, DOI 10.1016/j.aim.2006.07.010

\bibitem{ref7}{\it Y.C. Su, R.B. Zhang}:
 Cohomology of Lie superalgebras $\mathfrak{sl}_{m\mid n}$ and $\mathfrak{osp}_{2\mid2n}$. Proc. Lond. Math. Soc. {\bf94} (2007) 91-136. Zbl 1118.17005, MR2293466, DOI 10.1112/plms/pdl005

\bibitem{ref8}{\it B. Shu, C.W Zhang}:
 Restricted representations of the Witt superalgebras. J. Algebra. {\bf324} (2010), 652-672. Zbl 1217.17010, MR2651562, DOI 10.1016/J.JALGEBRA.2010.04.032

\bibitem{ref17}{\it B. Shu, C.W. Zhang}:
 Representations of the restricted Cartan type Lie superalgebra $W(m, n, 1)$. Algebr. Represent. Theor. {\bf14} (2011) 463-481. Zbl 1281.17020, MR2785918, DOI 10.1007/s10468-009-9198-6

\bibitem{ref9}{\it S.J. Wang, W.D. Liu}:
 The first cohomology of $\mathfrak{sl}_{2\mid1}$ with coefficients in $\chi$-reduced Kac modules and simple modules.
J. Pure Appl. Algebra. {\bf224} (2020), 106403. Zbl 1472.17066, MR4104490, DOI 10.1016/j.jpaa.2020.106403

\bibitem{ref13}{\it W.Q. Wang, L. Zhao}:
 Representations of Lie superalgebras in prime characteristic I. Proc. Lond. Math. Soc. {\bf99} (2009) 145-167. Zbl 1176.17013, MR2520353, DOI 10.1112/plms/pdn057

\bibitem{ref14}{\it W.Q. Wang, L. Zhao}:
 Representations of Lie Superalgebras in Prime Characteristic II: The queer series. J. Pure Appl. Algebra. {\bf215} (2011) 2515-2532. Zbl 1225.17021, MR2793955, DOI 10.1016/j.jpaa.2011.02.011

\bibitem{ref12}{\it L.S. Zheng}:
 On cohomology of modular Lie superalgebras. J. East China Norm. Univ. Natur. Sci. Ed. {\bf4} (2009) 82-91. Zbl 1212.17014, MR2555124

\bibitem{ref16}{\it C.W. Zhang}:
 On the simple modules for the restricted Lie superalgebra $\mathfrak{sl}(n|1)$. J. Pure Appl. Algebra. {\bf213} (2009) 756-765. Zbl 1230.17013, MR2494368, DOI 10.1016/j.jpaa.2008.09.005

 \bibitem{ref15}{\it L. Zhao}:
 Representations of Lie superalgebras in prime characteristic III. Pacific J. Math. {\bf248} (2010) 493-510. Zbl 1225.17022, MR2741259, DOI 10.2140/pjm.2010.248.493
\end{thebibliography}}
\end{document}